%% Notes on flat fronts in hyperbolic space
%% J Dubois, U Hertrich-Jeromin, G Szeswieczek
%%%%%%%%%%%%%%%%%%%%%%%%%%%%%%%%%%%%%%%%%%%%%%%%%%%%%%%%%%%%%%%

%%%%%%%%%%%%%%%%%%%%%%%%%%%%%%%%%%%%%%%%%%%%%%%%%%%%%%%%%%%%%%%
\def\ifudf#1{\expandafter\ifx\csname #1\endcsname\relax}
\newif\ifpdf \ifudf{pdfoutput}\pdffalse\else\pdftrue\fi
\ifpdf\pdfpagewidth=210mm \pdfpageheight=297mm
 \else\fi
\hsize=149mm \vsize=245mm \hoffset=5mm \voffset=0mm
\parskip=1ex minus .3ex \parindent=0pt
\everydisplay={\textstyle}
\font\tbb=bbmsl10 \font\sbb=bbmsl10 scaled 700
 \font\fbb=bbmsl10 scaled 500
\newfam\bbfam \textfont\bbfam=\tbb \scriptfont\bbfam=\sbb
\scriptscriptfont\bbfam=\fbb \def\bb{\fam\bbfam}%
\def\fn[#1]{\font\TmpFnt=#1\relax\TmpFnt\ignorespaces}
\def\em{\expandafter\ifx\the\font\tensl\rm\else\sl\fi}
\def\\{\hfill\break} \def\dt{1 Jun 2021}
\newtoks\secNo \secNo={0} \newcount\ssecNo \ssecNo=0
\def\section #1#2\par{\goodbreak\vskip 3ex\noindent
 \global\secNo={#1}\global\ssecNo=0 \global\eqnNo=0
 {\fn[cmbx10 scaled 1200]#1#2}\vglue 1ex}
\def\subsection #1.{\goodbreak\vskip 3ex\noindent
 \global\advance\ssecNo by 1 % \global\eqnNo=0
 {\fn[cmbx10]\the\secNo.\the\ssecNo.~#1}\vglue 1ex}
\def\href#1<#2>{\leavevmode
 \ifpdf\pdfstartlink attr {/Border [0 0 0 ]} goto name {#1}\fi
 {#2}\ifpdf\pdfendlink\fi}
\def\label@#1:#2@{\ifudf{#1}
 \expandafter\xdef\csname#1\endcsname{#2}\else
 \errmessage{label #1 already in use!}\fi}
\label@holom:3.1@ \label@W:3.2@ \label@XN:3.3@
\label@flatW:3.4@ \label@crH:3.5@ \label@ab:3.6@
\label@xy:3.7@ \label@crv:3.8@ \label@cc:3.9@
\label@mc11:3.10@ \label@mc12:3.11@ \label@bl29:1@
\label@bosu08:2@ \label@br87:3@ \label@bjr10:4@
\label@bjr12:5@ \label@bjr18:6@ \label@bcjpr20:7@
\label@ce92:8@ \label@cl12:9@ \label@gmm00:10@
\label@imdg:11@ \label@jeho17:12@ \label@jesz21:13@
\label@hrsy12:14@ \label@pe20:15@ \label@ppy21:16@
\label@roya16:17@
\newwrite\lbl \immediate\openout\lbl=\jobname.lbl
\def\@#1:#2@{\ifpdf\pdfdest name {#1} xyz\fi {#2}%
 \message{@#1:#2@}\immediate\write\lbl{\string\label@#1:#2@}}
\def\:#1:{\href#1<\ifudf{#1}??\else\csname#1\endcsname\fi>}
\newcount\refNo \refNo=0
\def\refitem#1 {\global\advance\refNo by 1
 \item{\@#1:\number\refNo@.}}
\let\pleqno\eqno \newcount\eqnNo \eqnNo=0
\def\eqno#1$${\global\advance\eqnNo by 1
 \pleqno{\rm(\@#1:\the\secNo.\number\eqnNo@)}$$}
\def\note #1.{\par\underbar{#1}.\enspace}
\def\C{{\bb C}} \def\CP{{\bb C}{\sl P}}
\def\P{{\bb P}} \def\R{{\bb R}} 
 
\def\Gl{{\rm Gl}} \def\Sl{{\rm Sl}} \def\Herm{{\rm Herm}}
\def\CR{{\rm cr}} \def\tr{\mathop{\rm tr}}
\def\<#1>{\langle #1\rangle}
%%%%%%%%%%%%%%%%%%%%%%%%%%%%%%%%%%%%%%%%%%%%%%%%%%%%%%%%%%%%%%%
%% top matter
%%%%%%%%%%%%%%%%%%%%%%%%%%%%%%%%%%%%%%%%%%%%%%%%%%%%%%%%%%%%%%%
\newtoks\title \newtoks\stitle \newtoks\author
\newtoks\status \newtoks\funding
\title={Notes on flat fronts in hyperbolic space}
\stitle={Notes on flat fronts}
\author={J Dubois, U Hertrich-Jeromin, G Szewieczek}
\status={arXiv version}
\ifpdf\pdfoutput=1\pdfadjustspacing=1\pdfinfo{%
 /Title (\the\title) /Author (\the\author) /Date (\dt)}\fi
\headline={\hfil} \footline={{\fn[cmr7]\hfil-\folio-\hfil}}
\centerline{{\fn[cmbx10 scaled 1440]\the\title}}\vglue .2ex
\centerline{{\fn[cmr7]\the\author}}\vglue 3em plus 3ex
\centerline{\vtop{\hsize=.8\hsize{\bf Abstract.}\enspace
 We give a short introduction to discrete flat fronts in
 hyperbolic space and prove that any discrete flat front
 in the mixed area sense admits a Weierstrass representation.
}}\vglue 2em
\centerline{\vtop{\hsize=.8\hsize{\bf MSC 2010.}\enspace
 {\it 53C42\/}, {\it 53A10\/}, 53A30, 37K25, 37K35
}}\vglue 1em
\centerline{\vtop{\hsize=.8\hsize{\bf Keywords.}\enspace
 flat front; hyperbolic space; linear Weingarten surface;
 discrete surface; curved flat; Darboux transformation;
 Weierstrass representation.
}}\vglue 3em
%%%%%%%%%%%%%%%%%%%%%%%%%%%%%%%%%%%%%%%%%%%%%%%%%%%%%%%%%%%%%%%

%%%%%%%%%%%%%%%%%%%%%%%%%%%%%%%%%%%%%%%%%%%%%%%%%%%%%%%%%%%%%%%
%% introduction
%%%%%%%%%%%%%%%%%%%%%%%%%%%%%%%%%%%%%%%%%%%%%%%%%%%%%%%%%%%%%%%
\section 1. Introduction
%%%%%%%%%%%%%%%%%%%%%%%%%%%%%%%%%%%%%%%%%%%%%%%%%%%%%%%%%%%%%%%

%%%%%%%%%%%%%%%%%%%%%%%%%%%%%%%%%%%%%%%%%%%%%%%%%%%%%%%%%%%%%%%
%% motivation: smooth
Flat fronts in hyperbolic space
--- that is, intrinsically flat surfaces in hyperbolic space
that may have certain types of singularities ---
yield an interesting class of surfaces/fronts for their
relations to various areas in differential geometry
and, in particular, also for the study of singularities
of special surfaces:
they admit a Weierstrass type representation [\:gmm00:],
that is of a rather geometric nature as it is directly linked
to the geometry of a pair of hyperbolic Gauss maps that form
a curved flat in the integrable systems sense [\:bjr10:].
For example, this observation establishes a relation with
Darboux pairs of minimal surfaces in Euclidean space
[\:jeho17:].
In a sphere geometric context, these hyperbolic Gauss maps
are related to a pair of isothermic sphere congruences that
put flat fronts in the context of Demoulin's \char 10-surfaces
and, in particular, of linear Weingarten surfaces [\:bjr12:];
on the other hand, flat fronts occur in parallel families
as orthogonal surfaces to a cyclic system spanned by the pair
of hyperbolic Gauss maps, cf [\:bjr10:].

%%%%%%%%%%%%%%%%%%%%%%%%%%%%%%%%%%%%%%%%%%%%%%%%%%%%%%%%%%%%%%%
%% motivation: discrete
Two approaches have been taken to discretize flat fronts in
hyperbolic space, in an integrable systems sense:
firstly, by discretizing their Weierstrass representation
 in [\:hrsy12:],
secondly, as particular linear Weingarten surfaces in
 a hyperbolic space form in [\:bjr18:].
A relation has been established by showing that the nets
obtained by the discrete Weierstrass-type representation
satisfy the second definition [\:roya16:],
though the converse seems not to be established yet.
However, in [\:jesz21:] we recently related discrete flat
fronts in the former sense to Darboux pairs
(curved flats of [\:cl12:])
of holomorphic maps into the $2$-sphere,
which is a building block towards proving the converse.

%%%%%%%%%%%%%%%%%%%%%%%%%%%%%%%%%%%%%%%%%%%%%%%%%%%%%%%%%%%%%%%
%% plan/outline
These notes grew out of the first author's MSc thesis,
supervised by the other authors:
we aim to give a relatively gentle and low-tech introduction
to the geometry of (discrete) flat fronts in hyperbolic space
and,
probably as our single new result,
to prove that Darboux pairs (curved flats) of discrete
holomorphic maps admit a Weierstrass-type representation.
As a consequence, any discrete linear Weingarten net
in hyperbolic space with mixed area Gauss curvature
$K\equiv 1$ then admits a Weierstrass-type representation.

%%%%%%%%%%%%%%%%%%%%%%%%%%%%%%%%%%%%%%%%%%%%%%%%%%%%%%%%%%%%%%%
%% acknowledgements
 {\it Acknowledgements.\/}
This work would not have been possible without the valuable
and enjoyable discussions around the subject area with
F Burstall and W Rossman.
Furthermore, we gratefully acknowledge financial support
from the Austrian Science Fund FWF through the
FWF/JSPS Joint Project I3809 ``Geometric shape generation''.

%%%%%%%%%%%%%%%%%%%%%%%%%%%%%%%%%%%%%%%%%%%%%%%%%%%%%%%%%%%%%%%
%% hyperbolic geometry
%%%%%%%%%%%%%%%%%%%%%%%%%%%%%%%%%%%%%%%%%%%%%%%%%%%%%%%%%%%%%%%
\section 2. Hyperbolic geometry
%%%%%%%%%%%%%%%%%%%%%%%%%%%%%%%%%%%%%%%%%%%%%%%%%%%%%%%%%%%%%%%

%%%%%%%%%%%%%%%%%%%%%%%%%%%%%%%%%%%%%%%%%%%%%%%%%%%%%%%%%%%%%%%
%% motivation
From a differential geometric point of view it is convenient
to embed hyperbolic space into a (flat) vector space geometry,
so that covariant differentiation becomes just ordinary
differentiation followed by a projection:
thus we consider hyperbolic space as one of the sheets
of a two-sheeted hyperboloid in a Minkowski space
$\R^{3,1}=(\R^4,(.,.))$:
$$
  H^3 = \{ x\in\R^{3,1} \,|\, (x,x)=-1,\, (x,o)<0 \},
$$
where $o\in H^3$ denotes some unit timelike vector,
that distinguishes one sheet of the hyperboloid.
To visualize objects in this $3$-dimensional hyperbolic
geometry we then project to a (non-flat) Poincar\'e model,
where $H^3$ is conformally embedded into Euclidean $3$-space,
for example, by means of the stereographic projection into
the unit ball of $\R^3=\{o\}^\perp\subset\R^{3,1}$:
$$
  \sigma:H^3\to\R^3,\enspace
  x \mapsto \sigma(x) := {x+o\,(o,x)\over 1-(o,x)}.
$$
In what follows we discuss two enhancements of this setup
that provide
for a better geometric handling of flat surfaces and fronts
 in hyperbolic space,
and
for an algebraic setting that caters for their Weierstrass
 type representation,
respectively.

%%%%%%%%%%%%%%%%%%%%%%%%%%%%%%%%%%%%%%%%%%%%%%%%%%%%%%%%%%%%%%%
%% model: Lie sphere geometry
%%%%%%%%%%%%%%%%%%%%%%%%%%%%%%%%%%%%%%%%%%%%%%%%%%%%%%%%%%%%%%%
\subsection Hyperbolic geometry as a subgeometry
 of Lie sphere geometry.
%%%%%%%%%%%%%%%%%%%%%%%%%%%%%%%%%%%%%%%%%%%%%%%%%%%%%%%%%%%%%%%
%% motivation/references
Hyperbolic geometry, as any other space form geometry,
may be considered as a subgeometry of Lie sphere geometry:
the classical incarnation of contact geometry.
This approach is, in particular, well suited to investigate
fronts, that is, surfaces $x:\Sigma^2\to H^3$ that may have
certain types of singularities but still admit a smooth
unit normal field (Gauss map) $n:\Sigma^2\to S^{2,1}$
such that the pair $(x,n)$ is an immersion.
Our short description here will be tailored to
the specific case of hyperbolic geometry,
for more general and comprehensive introductions
see [\:bl29:] or [\:ce92:].

%%%%%%%%%%%%%%%%%%%%%%%%%%%%%%%%%%%%%%%%%%%%%%%%%%%%%%%%%%%%%%%
%% general setup
Thus consider
$\R^{3,1}\subset\R^{4,2}=\R^{1,1}\oplus_\perp\R^{3,1}$
and let $(p,q)$ denote an orthonormal basis of $R^{1,1}$,
that is, $-(p,p)=(q,q)=1$ and $p\perp q$;
hence $\R^{3,1}=\<p,q>^\perp$,
where $\<\dots>$ denotes the linear span of vectors.
We now consider the affine quadrics
$$\matrix{
  Q^3 := \{ y\in{\cal L}\,|\,(y,p)=0,\,(y,q)=-1 \}, \cr
  P^3 := \{ y\in{\cal L}\,|\,(y,p)=-1,\,(y,q)=0 \}, \cr
   }\quad{\rm where}\enspace
  {\cal L} := \{ y\in\R^{4,2}\,|\,(y,y)=0 \}
$$
denotes the light cone of $\R^{4,2}$.
Clearly, $P^3$ and $Q^3$ are, as offsets of the unit $1$-
and $2$-sheeted hyperboloids in $\R^{3,1}$ by $p$ and $-q$,
respectively, isometric to these hyperboloids.
In particular, each component of $Q^3$ yields a model
of hyperbolic $3$-space in the light cone
${\cal L}\subset\R^{4,2}$.
The common infinity boundary of these two hyperbolic spaces
is formed by the null directions $\<x>$
in $\R^{3,1}\subset\R^{4,2}$,
where $x\in{\cal L}\cap\<p,q>^\perp$,
and inherits a $2$-dimensional M\"obius geometry from the
hyperbolic geometry of $Q^3$ resp the group of Lie sphere
transformations that fix $p$ and $q$.

%%%%%%%%%%%%%%%%%%%%%%%%%%%%%%%%%%%%%%%%%%%%%%%%%%%%%%%%%%%%%%%
%% spheres in hyperbolic space
More generally, points $\<s>$ of the Lie quadric
(the projective light cone of $\R^{4,2}$, i.e., $s\in{\cal L}$)
represent oriented spheres:
in hyperbolic geometry these come in various flavours;
the ``point spheres'', $s\perp p$, of either sheet of $Q^3$
or of their common infinity boundary have been discussed above.
Accordingly, we call $p$ the ``point sphere complex''.

A hyperbolic ``distance sphere''
with centre $c\in H^3\subset\R^{4,2}$ and radius $r\in\R$
is represented by $\<s>$ with $s=p\,\sinh r-q\,\cosh r+c$:
the intersection $Q^3\cap\<s>^\perp$ consists of points
in distance $r$ from $c-q$,
as evaluating $0=(x(t),s)$ on an arc-length parametrization
$t\mapsto x(t)=c\,\cosh t+c'\,\sinh t-q$
of a geodesic in $Q^3$ emanating from $c-q$ in a unit direction
$c'\in T_{c-q}Q^3=T_cH^3$ shows.
Note that a sign swap of the radius $r$ corresponds
to reflection of $s$ in $\<p>^\perp$:
this reflection yields the orientation reversal;
point spheres ($r=0$) are invariant under this reflection.
A reflection in $\<q>^\perp$, on the other hand,
swaps the sheets of $Q^3$, hence moving a distance sphere from
one hyperbolic space into the other.

Spheres represented by $\<p\,\cosh r+q\,\sinh r+n>$
with $n\in S^{2,1}$ and $r\in\R$ intersect the infinity sphere
of the two sheets of $Q^3$:
$\<n,p,q>^\perp\cap{\cal L}\neq\{0\}$.
In particular, if $r=0$ then $p+n\in P^3$ is invariant
under a reflection in $\<q>^\perp$,
that is, under an inversion in the infinity boundary,
hence represents a plane in both hyperbolic spaces.
Accordingly, we call $q$ the ``plane (sphere) complex''.

%%%%%%%%%%%%%%%%%%%%%%%%%%%%%%%%%%%%%%%%%%%%%%%%%%%%%%%%%%%%%%%
%% contact of spheres/contact elements
Given a point $x\in H^3$ and a unit vector $n\in T_xH^3$,
any sphere through $x$ with normal $n$ at $x$ is represented
by $s\in\<x-q,n+p>$,
for example, all distance spheres $(x-q)\cosh r+(n+p)\sinh r$
touch at the point $x$.
Note that $\<x-q,n+p>\subset\R^{4,2}$ is a null $2$-plane,
hence defines a line in the Lie quadric.
Such lines will be called ``contact elements''.

Two spheres $\<s>$ and $\<s'>$
``touch'' (are in oriented contact)
if they span a contact element, that is, if $(s,s')=0$.
In particular, a sphere is a horosphere if it touches either
of the oriented infinity spheres $\<q\pm p>$ of
the two hyperbolic sheets of $Q^3$:
thus a horosphere is represented by $h=q\pm p+x$,
where $x\in{\cal L}\cap\<p,q>^\perp$ yields
(homogeneous coordinates of)
the point of contact with the infinity boundary.

%%%%%%%%%%%%%%%%%%%%%%%%%%%%%%%%%%%%%%%%%%%%%%%%%%%%%%%%%%%%%%%
%% Lie sphere transformations
Lie sphere transformations are given by the orthogonal
transformations of $\R^{4,2}$, hence act on the space of
oriented spheres (the Lie quadric) and preserve contact.
Hyperbolic isometries are those Lie sphere transformations
that preserve the point and plane sphere complexes,
$\<p>$ and $\<q>$,
and that do not swap the sheets of $Q^3$.
Of interest will also be the ``parallel transformations'',
given by Lorentz boosts in $q$-direction,
$
  (p,q) \mapsto (p,q)\left(
   {\phantom{-}\cosh s\atop -\sinh s}\,
   {-\sinh s\atop\phantom{-}\cosh s}
  \right)
$,
that shift the radius of all distance spheres:
$$
  p\,\sinh r-q\,\cosh r+c
   \mapsto
  p\,\sinh(r+s)-q\,\cosh(r+s)+c.
$$

%%%%%%%%%%%%%%%%%%%%%%%%%%%%%%%%%%%%%%%%%%%%%%%%%%%%%%%%%%%%%%%
%% model: Pauli matrices/Hermitian forms
%%%%%%%%%%%%%%%%%%%%%%%%%%%%%%%%%%%%%%%%%%%%%%%%%%%%%%%%%%%%%%%
\subsection Hyperbolic geometry using Pauli matrices.
%%%%%%%%%%%%%%%%%%%%%%%%%%%%%%%%%%%%%%%%%%%%%%%%%%%%%%%%%%%%%%%
%% motivation/references
Weierstrass-type representations of surfaces in hyperbolic
space typically employ a model of hyperbolic geometry that
is based on Hermitian matrices (or forms),
using an explicit incarnation of the (universal) covering of
the group of hyperbolic motions (the restricted Lorentz group)
by the special linear group $\Sl(2,\C)$,
cf [\:br87:], [\:gmm00:] and [\:pe20:, Sect~4.2].
Here we will briefly recall this model of hyperbolic geometry,
for completeness.

%%%%%%%%%%%%%%%%%%%%%%%%%%%%%%%%%%%%%%%%%%%%%%%%%%%%%%%%%%%%%%%
%% Minkowski space-time via Hermitian matrices
Thus consider the (real) $4$-dimensional space of Hermitian
$2\times2$-matrices $\Herm(2)$:
eqipped with the quadratic form $-\det:\Herm(2)\to\R$
this space becomes a $4$-dimensional real Minkowski space,
$(\Herm(2),-\det) \cong \R^{3,1}$,
that $\Sl(2,\C)$ acts transitively on as its restricted Lorentz group via
$$
  \Sl(2,\C)\times\Herm(2)\to\Herm(2), \enspace
  (G,X)\mapsto G\cdot X=GXG^\ast.
$$
Together with the identity matrix the Pauli matrices form
an orthonormal basis $(E_0,E_1,E_2,E_3)$ of $\Herm(2)$,
$$
  E_0 = \left({1\atop 0}\,{0\atop\vphantom{-}1}\right),\enspace
  E_1 = \left({0\atop 1}\,{1\atop\vphantom{-}0}\right),\enspace
  E_2 = \left({0\atop i}\,{-i\atop\phantom{-}0}\right),\enspace
  E_3 = \left({1\atop 0}\,{\phantom{-}0\atop-1}\right).
$$
Hyperbolic geometry is then modelled on
$
  H^3 \cong \{ X\in\Herm(2)\,|\,\det X=1,\enspace \tr X>0\}
$
with the group $\Sl(2,\C)$ acting on $H^3$ by hyperbolic
motions;
note that orientation reversing hyperbolic isometries
(e.g., hyperbolic reflections)
are not obtained from the action of $\Sl(2,\C)$ on $\Herm(2)$.

%%%%%%%%%%%%%%%%%%%%%%%%%%%%%%%%%%%%%%%%%%%%%%%%%%%%%%%%%%%%%%%
%% Moebius infinity boundary geometry
In this setting of hyperbolic geometry, its boundary M\"obius
geometry comes in two flavours:
as the M\"obius geometry of restricted Lorentz transformations
 on the projective light cone of $\R^{3,1}$;
and
as the geometry of the complex projective line $\CP^1$,
 acted upon by complex M\"obius transformations.
The relation between the two models is established by
the Hermitian dyadic square,
$$
  \CP^1\ni \<x>\mapsto\<xx^\ast> \in\P({\cal L}^3),
   \enspace{\rm where}\enspace
  x \in\C^2
   \enspace{\rm and}\enspace
  {\cal L}^3=\{ X\in\Herm(2)\,|\,\det X=0 \},
$$
which is obviously compatible with the actions of $\Sl(2,\C)$
on $\CP^1$ by linear fractional transformations and
on the projective light cone $\P({\cal L}^3)$ of $\Herm(2)$
 by the restricted Lorentz group.
Also note that this map has an inverse as every singular
Hermitian matrix is a Hermitian dyadic square.

%%%%%%%%%%%%%%%%%%%%%%%%%%%%%%%%%%%%%%%%%%%%%%%%%%%%%%%%%%%%%%%
%% infinity endpoints of hyperbolic geodesics
For the hyperbolic Gauss maps of a surface or front in $H^3$
we will later be interested to determine the points,
where a geodesic emanating from a point $C\in H^3$ with a
unit initial velocity $C'\in T_CH^3$ hits the infinity sphere.
Using, as before, an arc-length parametrization we take the
limits
$$
  \lim_{t\to\pm\infty}\<C\,\cosh t+C'\,\sinh t> = \<C\pm C'>
$$
and observe that $-\det(C\pm C')=(C,C)\pm2(C,C')+(C',C')=0$,
that is, $C\pm C'\in{\cal L}^3$ indeed.
If, further, $C=G\cdot E_0$ and $C'=G\cdot E_3$ are given as
images of the basis vectors $E_0$ and $E_3$ under a hyperbolic
rotation $G=\left({a\atop b}\,{c\atop d}\right)$,
then
$$
  C+C' = G(E_0+E_3)G^\ast
  = 2\left({a\atop b}\right)\left({a\atop b}\right)^\ast
   \enspace{\rm and}\enspace
  C-C' = G(E_0-E_3)G^\ast
  = 2\left({c\atop d}\right)\left({c\atop d}\right)^\ast,
$$
that is, the columns of the matrix $G\in\Sl(2,\C)$ yield
(homogeneous coordinates for) the points at infinity of the
geodesic, as points in $\CP^1$.

%%%%%%%%%%%%%%%%%%%%%%%%%%%%%%%%%%%%%%%%%%%%%%%%%%%%%%%%%%%%%%%
%% discrete flat fronts
%%%%%%%%%%%%%%%%%%%%%%%%%%%%%%%%%%%%%%%%%%%%%%%%%%%%%%%%%%%%%%%
\section 3. Discrete flat fronts in hyperbolic space
%%%%%%%%%%%%%%%%%%%%%%%%%%%%%%%%%%%%%%%%%%%%%%%%%%%%%%%%%%%%%%%

%%%%%%%%%%%%%%%%%%%%%%%%%%%%%%%%%%%%%%%%%%%%%%%%%%%%%%%%%%%%%%%
%% motivation
Discrete flat fronts in hyperbolic space have been approached
from two rather different angles:
{\parindent=2em
\item{$\bullet$}
in [\:hrsy12:, Sect~4.3] they were defined in terms of
 a Weierstrass representation and shown to be principal
 (circular) nets (Thm 4.6 and Rem 4.7 in [\:hrsy12:]);
\item{$\bullet$}
in [\:bjr18:, Expl~4.3] they were investigated in terms
 of their mixed area (extrinsic) Gauss curvature,
 as a special instance of discrete linear Weingarten surfaces.
\par}
A relation was established in [\:roya16:, Sect~5],
where the discrete flat fronts of [\:hrsy12:] were shown
to be discrete linear Weingarten surfaces with a constant
mixed area Gauss curvature $K\equiv 1$,
see [\:roya16:, Prop~5.1 and Lemma 5.3].
Further, in [\:ppy21:], the Weierstrass-type representation
of [\:hrsy12:] is given a sphere geometric interpretation,
cf [\:ppy21:, Thm~4.9 and Expl~4.8].
Another sphere geometric approach to flat fronts,
in terms of (discrete) cyclic systems,
has been discussed in [\:jesz21:, Sect~4.2],
where a focus is on the geometry of the two hyperbolic
Gauss maps $h^\pm$ of a flat front
---
which draws a connection to the discrete curved flats
of [\:cl12:, Chap~8 (Thm~8.4.1)].

Here we aim to knit the various threads together,
to obtain further insight into the geometry of (discrete)
flat fronts and, in particular, to prove that every flat
front in the sense of [\:bjr18:] admits a Weierstrass-type
representation of [\:hrsy12:].
In the process we will present some short and simple proofs
for some of the results described above.

%%%%%%%%%%%%%%%%%%%%%%%%%%%%%%%%%%%%%%%%%%%%%%%%%%%%%%%%%%%%%%%
%% Weierstrass representation
%%%%%%%%%%%%%%%%%%%%%%%%%%%%%%%%%%%%%%%%%%%%%%%%%%%%%%%%%%%%%%%
\subsection A Weierstrass-type representation.
%%%%%%%%%%%%%%%%%%%%%%%%%%%%%%%%%%%%%%%%%%%%%%%%%%%%%%%%%%%%%%%
To set the scene we fix some notations,
cf [\:bjr18:] or [\:bcjpr20:]:
we shall consider maps that ``live'' on cells
of some dimension of a quadrilateral cell complex,
$\Sigma^2=(\Sigma^2_0,\Sigma^2_1,\Sigma^2_2)$, where
$\Sigma^2_0$ denotes the set of vertices,
$\Sigma^2_1$ the set of (oriented) edges, and
$\Sigma^2_2$ the set of (oriented) quadrilateral faces;
the elements of $\Sigma^2_1$ and $\Sigma^2_2$ will usually
be denoted by their vertices, i.e.,
$(ij)\in\Sigma^2_1$ denotes an edge from $i\in\Sigma^2_0$
 to $j\in\Sigma^2_0$, and
$(ijkl)\in \Sigma^2_2$ denotes an oriented face with vertices
 $i,j,k,l\in \Sigma^2_0$
 and edges $(ij),(jk),(kl),(li)\in \Sigma^2_1$.

We will generally assume $\Sigma^2$ to
be {\em simply connected\/},
that is,
any two points $i,j\in\Sigma^2_0$ can be connected by
an (edge-)path and any closed (edge-)path is null-homotopic,
via ``face-flips'' and
via dropping single edge ``return trips'',
cf [\:bcjpr20:, Sect~2.3].

%%%%%%%%%%%%%%%%%%%%%%%%%%%%%%%%%%%%%%%%%%%%%%%%%%%%%%%%%%%%%%%
%% derived edge function and derivative/edge-labelling
Given a vector-valued function $i\mapsto g_i$ on vertices,
its {\em derived function\/} and discrete {\em derivative\/}
will be denoted by
$$
  g_{ij} := {g_i+g_j\over 2}
   \enspace{\rm and}\enspace
  dg_{ij} := g_j-g_i,
$$
so that a Leibniz rule $d(fg)_{ij}=df_{ij}g_{ij}+f_{ij}dg_{ij}$
(for any product of vectors) holds;
note that
the derived edge function is a function on unoriented edges
while
the derivative is a discrete $1$-form,
$$
  g_{ij} = g_{ji}
   \enspace{\rm while}\enspace
  dg_{ij} = -dg_{ji}.
$$

%%%%%%%%%%%%%%%%%%%%%%%%%%%%%%%%%%%%%%%%%%%%%%%%%%%%%%%%%%%%%%%
%% discrete holomorphic functions/edge labelling
If $g:\Sigma^2_0\to\C$ is a complex-valued vertex function,
then we may compute the cross ratio on (quadrilateral)
faces,
$$
  \CR(g_i,g_j,g_k,g_l) :
   = {dg_{ij}\over dg_{jk}}{dg_{kl}\over dg_{li}};
$$
note that this does not define a map on faces,
but it does after factoring out the anharmonic group.
%% cf [\:imdg:, \S4.9.11].
Then $g$ will be called {\em holomorphic\/}
if it is a (totally umbilic) discrete isothermic net:
$$
  \forall(ijkl)\in\Sigma^2_2:
  \CR(g_i,g_j,g_k,g_l) = {a_{ij}\over a_{jk}},
\eqno holom$$
where $a:\Sigma^2_1\to\R$ denotes
a (real) {\em edge-labelling\/},
that is, a function on edges so that
$$
  \forall(ij)\in\Sigma^2_1:a_{ij}=a_{ji}
   \enspace{\rm and}\enspace
  \forall(ijkl)\in\Sigma^2_2:a_{ij}=a_{kl}.
$$
For regularity, we further assume that
$
  \forall(ijkl)\in\Sigma^2_2:
  \CR(g_i,g_j,g_k,g_l) \neq 0,1,\infty.
$

%%%%%%%%%%%%%%%%%%%%%%%%%%%%%%%%%%%%%%%%%%%%%%%%%%%%%%%%%%%%%%%
%% Weierstrass representation
We are now in a position to recall the representation
of discrete flat fronts from [\:hrsy12:, Sect 4.3]:

\proclaim Weierstrass representation.
Let $g:\Sigma^2_0\to\C$ be a discrete holomorphic function,
with edge-labelling $a:\Sigma^2_1\to\R$,
and fix
$t\in\R\setminus\{0,{1\over a_{ij}}\,|\,(ij)\in\Sigma^2_1\}$;
then there is a frame $F:\Sigma^2_0\to\Sl(2,\C)$, satisfying
$$
  F_j = F_iW_{ij}
   \enspace{\rm with}\enspace
  W_{ij} := \pmatrix{
   1 & dg_{ij} \cr
   {ta_{ij}\over dg_{ij}} & 1 \cr
   }\,{1\over\sqrt{1-ta_{ij}}},
\eqno W$$
that yields a $1$-parameter family of parallel discrete
flat fronts $(X,N):\Sigma^2_0\to H^3\times S^{2,1}$,
where
$$
  X := F\,(E_0\cosh s+E_3\sinh s)\,F^\ast
   \enspace{\sl and}\enspace
  N := F\,(E_0\sinh s+E_3\cosh s)\,F^\ast
\eqno XN$$
for $s\in\R$ and 
$
  E_0 = \left({1\atop 0}\,{0\atop\vphantom{-}1}\right),
  E_3 = \left({1\atop 0}\,{\phantom{-}0\atop-1}\right)
  \in\Herm(2),
$
as above.

%%%%%%%%%%%%%%%%%%%%%%%%%%%%%%%%%%%%%%%%%%%%%%%%%%%%%%%%%%%%%%%
%% Christoffel dual holomorphic net/Weierstrass connection
%% ...note that we use the conjugate of imdg:\S5.7.7
\note Remark.
By Christoffel's formula
$dg^\ast_{ij}={a_{ij}\over d\bar g_{ij}}$
with the (isothermic) dual
$g^\ast$ of $g$, cf [\:imdg:, \S5.7.7];
hence
$$
  W = \left({1\atop td\bar g^\ast}\,{dg\atop 1}\right)
      {1\over\sqrt{1-ta}}.
$$
Note that $W:\Sigma^2_1\to\Sl(2,\C)$ defines
a {\em discrete connection\/} on the vector bundle
$\Sigma^2_0\times\Herm(2)$, since
$$
  W_{ij}W_{ji}
  = \left({1\atop td\bar g^\ast_{ij}}\,{dg_{ij}\atop 1}\right)
    \left({1\atop-td\bar g^\ast_{ij}}\,{-dg_{ij}\atop 1}\right)
    {1\over\sqrt{1-ta_{ij}}^2}
  = \left({1\atop 0}\,{0\atop\vphantom{-}1}\right).
$$

%%%%%%%%%%%%%%%%%%%%%%%%%%%%%%%%%%%%%%%%%%%%%%%%%%%%%%%%%%%%%%%
%% proof: existence of frame
\note Integrability.
To substantiate the claims in relation to the above
Weierstrass-type representation, first note that
the frame $F$ indeed exists:
as $\Sigma^2$ is assumed to be simply connected,
the integrability of (\:W:) is equivalent to
the (metric) connection $W$ being flat
(cf [\:bcjpr20:, Prop~2.17]),
that is,
$$
  \forall(ijkl)\in\Sigma^2_2:
  W_{ij}W_{jk} = W_{il}W_{lk}.
\eqno flatW$$
Using that $a:\Sigma^2_1\to\R$ is an edge-labelling,
the flatness (\:flatW:) of the connection $W$ reduces to
the integrability of the Christoffel equation for $g^\ast$,
and to
the cross ratio condition (\:holom:) for $g$ to be holomorphic,
cf [\:hrsy12:, Thm~4.6 (proof)] and [\:roya16:, Sect 4.3].

%%%%%%%%%%%%%%%%%%%%%%%%%%%%%%%%%%%%%%%%%%%%%%%%%%%%%%%%%%%%%%%
%% geometric propagation of Legendre map by reflections
\note Remark.
The Weingarten equation (\:W:) yields a propagation
for the pair $(X,N)$ via reflections.
Namely, let
$$
  Y_{ij} := \left(
    {|dg_{ij}|^2 \atop d\bar g_{ij}}\,
    { dg_{ij}    \atop      ta_{ij}}
   \right)
   {1\over |dg_{ij}|\,\sqrt{1-ta_{ij}}}
  \in\Herm(2)
$$
and observe that
$
  F_iY_{ij}F_i^\ast + F_jY_{ji}F_j^\ast
  = F_j(W_{ji}Y_{ij}+Y_{ji}W_{ij}^\ast)F_i^\ast
  = 0;
$
hence $R_{ij}:=F_iY_{ij}F_i^\ast$ yields a well defined
reflection on $\R^{3,1}\cong(\Herm(2),-\det)$,
$$
  X \mapsto \varrho_{ji}(X) :
  = X - 2R_{ij}(R_{ij},X)
  = X + R_{ij}{\det(X+R_{ij})-\det(X-R_{ij})\over 2}.
$$
As $E_k-2Y_{ij}(Y_{ij},E_k)=W_{ij}E_kW_{ij}^\ast$ for $k=0,3$
we find that these reflections transport the points between
adjacent geodesics from $X$ in direction $N$,
$$
  X_j = \varrho_{ji}(X_i)
   \enspace{\rm and}\enspace
  N_j = \varrho_{ji}(N_i)
   \enspace{\rm for}\enspace
  s\in\R.
$$
Note that, in accordance with [\:jesz21:, Sect~4.2],
$\<R_{ij}>=\<d(X\pm N)_{ij}>$.

%%%%%%%%%%%%%%%%%%%%%%%%%%%%%%%%%%%%%%%%%%%%%%%%%%%%%%%%%%%%%%%
%% Legendre map/circularity
\note Circularity.
In fact, the map $(X,N):\Sigma^2_0\to H^3\times S^{2,1}$
satisfies a discrete {\em Rodrigues equation\/} on edges,
$$
  0 = dN_{ij} + k_{ij}dX_{ij},
   \enspace{\rm where}\enspace
  k_{ij} = -{(dN_{ij},R_{ij})\over(dX_{ij},R_{ij})}
$$
is the {\em principal curvature\/} of the edge;
note that $k_{ij}=k_{ji}$.
To obtain a Lie geometric interpretation we embed
$\R^{3,1}\subset\R^{4,2}=\R^{3,1}\oplus_\perp\R^{1,1}$
with $\R^{1,1}=\<p,q>$, as before,
then consider the discrete {\em sphere congruences\/}
$\<x>$ and $\<n>$ given by the light cone maps
$$
%%  x := X-q, n := N+p: \Sigma^2_0\to{\cal L}\subset\R^{4,2}.
  X-q=:x: \Sigma^2_0\to Q^3
   \enspace{\rm and}\enspace
  N+p=:n: \Sigma^2_0\to P^3.
$$
Clearly, $dn=dN$ and $dx=dX$ satisfy the same Rodrigues
equation, hence adjacent contact elements intersect in
a common (edge-){\em curvature sphere\/}
$
  \<\kappa_{ij}> = \<x,n>_i\cap\<x,n>_j
$
--- which characterizes the contact element map
$\<x,n>$ as a discrete {\em Legendre map\/},
cf [\:bosu08:, Sect~3.5] or [\:bjr18:, Def~2.1].
As a consequence,
each map $X:\Sigma^2_0\to H^3$ resp $x:\Sigma^2_0\to Q^3$
of the parallel family is a {\em principal\/}
(or, {\em circular\/}) {\em net\/},
cf [\:hrsy12:, Thm~4.6] and [\:roya16:, Lemma~5.1].

%%%%%%%%%%%%%%%%%%%%%%%%%%%%%%%%%%%%%%%%%%%%%%%%%%%%%%%%%%%%%%%
%% parallel family via Lie sphere geometry
\note Remark.
As $x:\Sigma^2_0\to Q^3$ and $n:\Sigma^2_0\to P^3$, these are
the point sphere and tangent plane congruences of the Legendre
map $\<x,n>$, respectively.
Note how the aforementioned hyperbolic parallel transformations
of Sect~2.1 move through the parallel family in the Weierstrass
representation:
$$\matrix{
    \<X-(p\sinh s+q\cosh s),N+(p\cosh s+q\sinh s)> \hfill\cr
  = \<(X\cosh s+N\sinh s)-q,(X\sinh s+N\cosh s)+p>. \cr
}$$

%%%%%%%%%%%%%%%%%%%%%%%%%%%%%%%%%%%%%%%%%%%%%%%%%%%%%%%%%%%%%%%
%% flat front Gauss curvature
\note Curvature.
Finally, we argue that $(X,N):\Sigma^2_0\to H^3\times S^{2,1}$
defines a flat front;
to this end we compute face diagonals
 $\delta H^\pm_{ik}:=H^\pm_k-H^\pm_i$
for its hyperbolic Gauss maps,
 $H^\pm=X\pm N$,
then determine
its (extrinsic) mixed area {\em Gauss curvature\/} $K$
on faces from the equation (cf [\:bjr18:, Def~2.3])
$$
  A(H^+,H^-) = A(X,X) - A(N,N) = A(X,X)(1-K),
$$
where $A:\Sigma^2_2\to\Lambda^2\R^{3,1}$ denotes
the (vector-valued) mixed area,
a symmetric bilinear form given by its quadratic form
$$
  A(H,H)_{ijkl} :
  = {1\over 2}\delta H_{ik}\wedge\delta H_{jl}.
$$
Thus let $E^\pm:=E_0\pm E_3$, so that $H^\pm=FE^\pm F^\ast$,
and let $r:\Sigma^2_0\to\R$ be a function that factorizes
${|dg|^2\over a}:\Sigma^2_1\to\R$, so that 
$
  dg^\ast_{ij}
  = {a_{ij}\over d\bar g_{ij}}
  = {dg_{ij}\over r_ir_j}
$
(cf [\:bosu08:, Thm~2.31] or [\:bjr18:, (3.1)]);
hence compute
%$$\matrix{
%  W_{ij}W_{jk}E_+ &=& \pmatrix{
%    1+tdg_{ij}d\bar g^\ast_{jk} & 0 \cr
%    t(d\bar g^\ast_{ij}+d\bar g^\ast_{jk}) & 0 \cr
%   }{2\over\sqrt{~}\sqrt{~}}\hfill\cr
%  W_{ij}W_{jk}E_- &=& \pmatrix{
%    0 & dg_{ij}+dg_{jk} \cr
%    0 & 1+td\bar g^\ast_{ij}dg_{jk} \cr
%   }{2\over\sqrt{~}\sqrt{~}}\hfill\cr
%  W_{ij}E_+W_{kj}^\ast &=& \pmatrix{
%    1 & -tdg^\ast_{jk} \cr
%    td\bar g^\ast_{ij} & -t^2d\bar g^\ast_{ij}dg^\ast_{jk} \cr
%   }{2\over\sqrt{~}\sqrt{~}}\hfill\cr
%  W_{ij}E_-W_{kj}^\ast &=& \pmatrix{
%    -dg_{ij}d\bar g_{jk} & dg_{ij} \cr
%    -d\bar g_{jk} & 1 \cr
%   }{2\over\sqrt{~}\sqrt{~}}\hfill\cr
%}$$
$$\matrix{
  F_i^{-1}\delta H^+_{ik}(F^\ast_k)^{-1}
  &=& W_{ij}W_{jk}E^+ - (W_{kl}W_{li}E^+)^\ast \hfill\cr
%  &=& \pmatrix{
%     dg_{ij}d\bar g^\ast_{jk}-d\bar g_{kl}dg^\ast_{li} &
%    -(dg^\ast_{kl}+dg^\ast_{li}) \cr
%     d\bar g^\ast_{ij}+d\bar g^\ast_{jk} & 0 \cr
%   }{2t\over\sqrt{~}\sqrt{~}}\hfill\cr
  &=& \pmatrix{
   ({dg_{ij}\over dg_{jk}}{a_{jk}\over r_k}
   -{dg_{il}\over dg_{lk}}{a_{ij}\over r_i})\,r_k &
    \delta g^\ast_{ik} \cr
   \delta\bar g^\ast_{ik} & 0 \cr
   }{2t\over\sqrt{1-ta_{ij}}\sqrt{1-ta_{jk}}}, \hfill\cr
  F_i^{-1}\delta H^-_{jl}(F^\ast_k)^{-1}
  &=& W_{il}E^-W_{kl}^\ast - W_{ij}E^-W_{kj}^\ast \hfill\cr
%  &=& \pmatrix{
%    -dg_{il}d\bar g_{lk}+dg_{ij}d\bar g_{jk} &
%     dg_{il}-dg_{ij} \cr
%    -d\bar g_{lk} + d\bar g_{jk} & 0 \cr
%   }{2\over\sqrt{~}\sqrt{~}}\hfill\cr
  &=& \pmatrix{
   ({dg_{ij}\over dg_{jk}}a_{jk}r_j
   -{dg_{il}\over dg_{lk}}a_{ij}r_l)\,r_k & \delta g_{jl} \cr
   \delta\bar g_{jl} & 0 \cr
   }{2\over\sqrt{1-ta_{ij}}\sqrt{1-ta_{jk}}}. \hfill\cr
}$$
Since
$
  0 = {\delta g^\ast_{ik}\over\delta(1/r)_{ik}}
    + {\delta g_{jl}\over\delta r_{jl}},
$
from the Christoffel equation
(cf [\:bosu08:, (2.37)] or [\:bjr18:, (3.3)]),
we conclude that
$$
  \forall(ijkl)\in\Sigma^2_2:
  \delta H^+_{ik}\parallel\delta H^-_{jl}
   \enspace\Rightarrow\enspace
  A(H^+,H^-) \equiv 0
   \enspace\Rightarrow\enspace
  K \equiv 1,
$$
that is, $(X,N)$ has the extrinsic Gauss curvature of
a (smooth) flat surface in hyperbolic space,
cf [\:roya16:, Lemma~5.3].
Furthermore, as $\delta(x\pm n)\parallel\delta H^\pm$,
we learn that the lifts $x\pm n:\Sigma^2_0\to{\cal L}$
of the enveloped horosphere congruences $\<x\pm n>$ are
K\"onigs dual lifts of isothermic sphere congruences,
showing that the Legendre map $\<x,n>=\<x+n,x-n>$
is an \char 10-net in the sense of [\:bjr18:, Def~3.1]
 or [\:bcjpr20:, Def~6.1],
cf [\:roya16:, Prop~5.1].

%%%%%%%%%%%%%%%%%%%%%%%%%%%%%%%%%%%%%%%%%%%%%%%%%%%%%%%%%%%%%%%
%% summary
\proclaim Summary.
Any discrete flat front $(X,N):\Sigma^2_0\to H^3\times S^{2,1}$
obtained from the Weierstrass representation (\:XN:) is
a discrete linear Weingarten net with constant Gauss curvature
$K\equiv 1$.

%%%%%%%%%%%%%%%%%%%%%%%%%%%%%%%%%%%%%%%%%%%%%%%%%%%%%%%%%%%%%%%
%% Darboux pairs of holomorphic maps
%%%%%%%%%%%%%%%%%%%%%%%%%%%%%%%%%%%%%%%%%%%%%%%%%%%%%%%%%%%%%%%
\subsection Darboux pairs of holomorphic maps.
%%%%%%%%%%%%%%%%%%%%%%%%%%%%%%%%%%%%%%%%%%%%%%%%%%%%%%%%%%%%%%%
%% motivation
This last observation,
about the enveloped horosphere congruences $\<x\pm n>$
of a flat front $(x,n)$ in hyperbolic space being isothermic
with K\"onigs dual lifts,
leads to a relation between flat fronts and Darboux pairs
of discrete holomorphic maps, see [\:jesz21:, Thm~31]:
 {\sl the hyperbolic Gauss maps of a discrete flat front
 (in the mixed area sense)
 form a Darboux pair of holomorphic maps in the $2$-sphere;
 conversely,
 any Darboux pair of discrete holomorphic maps in a $2$-sphere
 yields a hyperbolic line congruence so that its orthogonal nets
 have mixed area Gauss curvature $K\equiv 1$\/}.

%%%%%%%%%%%%%%%%%%%%%%%%%%%%%%%%%%%%%%%%%%%%%%%%%%%%%%%%%%%%%%%
%% Darboux pair from flat front/Weierstrass representation
Clearly the first assertion of [\:jesz21:, Thm~31],
for discrete flat fronts in the sense of [\:bjr18:],
now follows directly also for discrete flat fronts
in the sense of [\:hrsy12:]:

\proclaim Lemma.
The two hyperbolic Gauss maps $\<H^\pm>$ of a
discrete flat front $(X,N):\Sigma^2_0\to H^3\times S^{2,1}$,
obtained from the Weierstrass representation (\:XN:),
form a Darboux pair in $S^2\cong\P({\cal L}^3)$.

Alternatively, this Lemma is more directly verified by
straightforward cross ratio computations:
using the $\CP^1$-versions $\<h^\pm>:\Sigma^2_0\to\CP^1$,
$h^\pm=Fe^\pm$ with
 $e^+=\left({1\atop0}\right)$ and
 $e^-=\left({0\atop1}\right)$,
of the hyperbolic Gauss maps,
the Weierstrass representation (\:W:) yields
$$\matrix{
  \CR(\<h^\pm_i>,\<h^\pm_j>,\<h^\pm_k>,\<h^\pm_l>)
%  &=& {\det(F_ie^\pm,F_je^\pm)\over\det(F_je^\pm,F_ke^\pm)}
%      {\det(F_ke^\pm,F_le^\pm)\over\det(F_le^\pm,F_ie^\pm)}
%     \hfill\cr
%  &=& {\det(F_ie^\pm,F_iW_{ij}e^\pm)\over
%       \det(F_iW_{ij}e^\pm,F_iW_{ij}W_{jk}e^\pm)}
%      {\det(F_iW_{il}W_{lk}e^\pm,F_iW_{il}e^\pm)\over
%       \det(F_iW_{il}e^\pm,F_ie^\pm)}
%     \hfill\cr
  &=& {\det(e^\pm,W_{ij}e^\pm)\over\det(e^\pm,W_{jk}e^\pm)}
      {\det(W_{lk}e^\pm,e^\pm)\over\det(W_{il}e^\pm,e^\pm)}
%     \hfill\cr
  &=& {a_{ij}\over 1-ta_{ij}}{1-ta_{jk}\over a_{jk}},
  \hfill\cr
  \CR(\<h^+_i>,\<h^+_j>,\<h^-_j>,\<h^-_i>)
%  &=& {\det(F_ie^+,F_je^+)\over\det(F_je^+,F_je^-)}
%      {\det(F_je^-,F_ie^-)\over\det(F_ie^-,F_ie^+)}
%     \hfill\cr
  &=& {\det(e^+,W_{ij}e^+)\over\det(e^+,e^-)}
      {\det(W_{ij}e^-,e^-)\over\det(e^-,e^+)}
%     \hfill\cr
  &=& -{ta_{ij}\over 1-ta_{ij}}.
  \hfill\cr
}\eqno crH$$
Thus $\<h^\pm>:\Sigma^2_0\to\CP^1\cong S^2$ form a Darboux pair
of isothermic nets, with cross ratio factorizing edge-labelling
$$
  b := -{a\over 1-ta}
   \enspace\Leftrightarrow\enspace
  1 = (1-ta)(1-tb).
\eqno ab$$

%%%%%%%%%%%%%%%%%%%%%%%%%%%%%%%%%%%%%%%%%%%%%%%%%%%%%%%%%%%%%%%
%% curved flat of holomorphic maps
\note Remark.
By [\:cl12:, Thm~8.4.1] we just confirmed that the point pair
map $(\<h^+>,\<h^->):\Sigma^2_0\to S^2\times S^2$ qualifies
as a discrete curved flat.
Observe how this resonates with the change of the cross ratio
factorizing edge-labelling (\:ab:),
cf [\:imdg:, \S5.7.34 and \S5.7.16],
and the fact that (\:W:) is a discretization of the smooth
curved flat system $dF=F\Phi$ with
$\Phi=\left({0\atop td\bar g^\ast}\,{dg\atop 0}\right)$,
cf [\:imdg:, \S5.5.20] and [\:bjr10:, Sect 3].

%%%%%%%%%%%%%%%%%%%%%%%%%%%%%%%%%%%%%%%%%%%%%%%%%%%%%%%%%%%%%%%
%% mission: Weierstrass representation from Darboux pair
\note Mission.
We would like to obtain the converse:
given
%% a flat front $(X,N)$ with hyperbolic Gauss maps $H^\pm$
%% or, equivalently,
a Darboux pair (curved flat) of discrete holomorphic maps
$\<h^\pm>:\Sigma^2_0\to\CP^1$ into a $2$-sphere $S^2\cong\CP^1$,
we would like to show the existence of
--- and derive a defining equation for ---
lifts $h^\pm$ so that $F=(h^+,h^-)$ satisfies (\:W:).

%%%%%%%%%%%%%%%%%%%%%%%%%%%%%%%%%%%%%%%%%%%%%%%%%%%%%%%%%%%%%%%
%% setup to accomplish mission
\note Setting.
Thus let $h^\pm:\Sigma^2_0\to\C^2$ denote (any lifts of)
the maps of a Darboux pair of discrete holomorphic maps
into the $2$-sphere $S^2\cong\CP^1$,
with parameter $t\in\R\setminus\{0\}$
and an edge-labelling $b:\Sigma^2_1\to\R$,
that is (cf [\:imdg:, \S5.7.12]),
$$\matrix{
  \forall(ijkl)\in\Sigma^2_2:&
   \CR(\<h^\pm_i>,\<h^\pm_j>,\<h^\pm_k>,\<h^\pm_l>)
    = {b_{ij}\over b_{jk}}
   \enspace{\rm and}\cr\hfill
  \forall(ij)\in\Sigma^2_1:&
   \CR(\<h^+_i>,\<h^+_j>,\<h^-_j>,\<h^-_i>)
    = tb_{ij}. \hfill\cr
}$$
Further, let $F':=(h^+,h^-):\Sigma^2_0\to\Gl(2,\C)$ and
$
  W':\Sigma^2_1\to\Gl(2,\C), \enspace
  (ij)\mapsto W'_{ij} := (F'_i)^{-1}F'_j.
$

%%%%%%%%%%%%%%%%%%%%%%%%%%%%%%%%%%%%%%%%%%%%%%%%%%%%%%%%%%%%%%%
%% cr computations
\note Cross ratios.
Writing $W'=\left({u\atop x}\,{v\atop y}\right)$,
the assumption $F':\Sigma^2_0\to\Sl(2,\C)$ and
the condition for the cross ratio of corresponding edges
yield, analogously to (\:crH:),
$$
  \left.\matrix{
   uy-vx = 1 \hfill\cr
   \phantom{uy}-vx = tb \cr
  }\right\}\enspace\Rightarrow\enspace\left\{\matrix{
   x = -{tb\over v}, \hfill\cr
   y = {1-tb\over u}; \cr
  }\right.
\eqno xy$$
hence the factorizing cross ratios of the isothermic nets
$h^\pm$ on $(ijkl)\in\Sigma^2_2$ reduce to a single equation,
$$
  {b_{ij}\over b_{jk}}
  = \CR(\<h^\pm_i>,\<h^\pm_j>,\<h^\pm_k>,\<h^\pm_l>)
  = {v_{ij}\over v_{jk}}{v_{lk}\over v_{il}}
   \enspace\Leftrightarrow\enspace
  b_{jk}{v_{ij}\over v_{jk}} = b_{ij}{v_{il}\over v_{lk}}.
\eqno crv$$

%%%%%%%%%%%%%%%%%%%%%%%%%%%%%%%%%%%%%%%%%%%%%%%%%%%%%%%%%%%%%%%
%% connection/mc computations
\note Maurer-Cartan equation.
Next we wheel out the conditions on $W'$ to define
a flat connection:
firstly, using (\:xy:), the condition $W'_{ji}=(W'_{ij})^{-1}$
on $W'$ to define a connection produces two equations
$$
  \left.\matrix{
   1 &=& u_{ij}u_{ji} - tb_{ij}{v_{ij}\over v_{ji}} \hfill\cr
   0 &=& u_{ij}u_{ji} + (1-tb_{ij}){v_{ij}\over v_{ji}} \cr
  }\right\}
   \enspace\Leftrightarrow\enspace
  \left\{\matrix{\hfill
   v_{ji} &=& -v_{ij}, \hfill\cr
   u_{ij}u_{ji} &=& 1-tb_{ij}, \cr
  }\right.
\eqno cc$$
in particular, we conclude that $v$ is a discrete $1$-form;
secondly, we use (\:xy:) and (\:crv:) to extract two equations
from the flatness $W'_{ij}W'_{jk}=W'_{il}W'_{lk}$ of $W'$,
$$
   u_{ij}u_{jk} - u_{il}u_{lk}
   = t(b_{jk}{v_{ij}\over v_{jk}} - b_{ij}{v_{il}\over v_{lk}})
   = 0
\eqno mc11$$
and
$$
    v_{ij}{1-tb_{jk}\over u_{jk}} + v_{jk}u_{ij}
  = v_{il}{1-tb_{ij}\over u_{lk}} + v_{lk}u_{il}.
\eqno mc12$$

%%%%%%%%%%%%%%%%%%%%%%%%%%%%%%%%%%%%%%%%%%%%%%%%%%%%%%%%%%%%%%%
%% gauge transformation
\note Gauge transformation.
Finally consider the remaining freedom $F=F'G=(h^+w,h^-w^{-1})$
for the choice of lifts $(h^+,h^-)$ of the Darboux pair
and the corresponding gauge transformation for
the induced connection,
$$
  W_{ij}
%  = F_i^{-1}F_j
  = G^{-1}_iW'_{ij}G_j
%  = \pmatrix{
%     u_{ij}{w_j\over w_i} & v_{ij}{1\over w_iw_j} \cr
%     x_{ij}w_iw_j & y_{ij}{w_i\over w_j} \cr
%    }
  = \pmatrix{
     {u_{ij}w_j\over w_i} &
     {v_{ij}\over w_iw_j} \cr
     -tb_{ij}{w_iw_j\over v_{ij}} &
     (1-tb_{ij}){w_i\over u_{ij}w_j} \cr
    },
$$
seeking lifts so that $W$ has the appropriate shape of (\:W:):
using (\:ab:) this leads to the difference equation
$$
  \left.\matrix{\hfill
    {u_{ij}w_j\over w_i} \cr
    (1-tb_{ij}){w_i\over u_{ij}w_j} \cr
    }\right\}
  = {1\over\sqrt{1-ta_{ij}}}
  = \sqrt{1-tb_{ij}}
   \enspace\Leftrightarrow\enspace
  w_j = {\sqrt{1-tb_{ij}}\over u_{ij}}\,w_i,
$$
the integrability of which is granted by (\:mc11:) and
the fact that $b$ is an edge-labelling.
Further, using the difference equation for $w$,
(\:cc:) and (\:mc12:) yield the existence of a potential
$g:\Sigma^2_0\to\C$ of the $1$-form
$$
  \Sigma^2_1 \ni (ij)\mapsto
   {v_{ij}\over w_iw_j\sqrt{1-tb_{ij}}} = dg_{ij} \in\C,
$$
% ...integrability:
% $$
%     {v_{ij}\over w_iw_j\sqrt{1-tb_{ij}}}
%   + {v_{jk}\over w_jw_k\sqrt{1-tb_{jk}}}
%   = {1\over w_iw_k\sqrt{1-tb_{ij}}\sqrt{1-tb_{jk}}}
%     \{ v_{ij}{1-tb_{jk}\over u_{jk}} + v_{jk}u_{ij} \}
% $$
which is holomorphic by (\:crv:),
with cross ratio factorizing edge-labelling $a=-{b\over 1-tb}$
by (\:ab:).
% ...finally, $d\bar g^\ast$:
% $$
%   -b_{ij}{w_iw_j\over v_{ij}}
%   = -{b_{ij}\over 1-b_{ij}}
%     {w_iw_j\sqrt{1-tb_{ij}}\over v_{ij}}
%     \sqrt{1-tb_{ij}}
%   = {a_{ij}\over dg_{ij}}
%     {1\over\sqrt{1-ta_{ij}}}
% $$

%%%%%%%%%%%%%%%%%%%%%%%%%%%%%%%%%%%%%%%%%%%%%%%%%%%%%%%%%%%%%%%
%% Weierstrass representation from Darboux pair
Thus we have proved:

\proclaim Lemma.
Any Darboux pair of discrete holomorphic maps in the $2$-sphere
$S^2\cong\CP^1$ admits a Weierstrass representation (\:W:),
that is, there exist lifts $h^\pm:\Sigma^2_0\to\C^2$ of the
two holomorphic maps $\<h^\pm>:\Sigma^2_0\to\CP^1$
so that $F:=(h^+,h^-):\Sigma^2_0\to\Sl(2,\C)$ satisfies (\:W:).

%%%%%%%%%%%%%%%%%%%%%%%%%%%%%%%%%%%%%%%%%%%%%%%%%%%%%%%%%%%%%%%
%% existence of Weierstrass representation for flat fronts
And, as the hyperbolic Gauss maps of a discrete flat front
in the mixed area sense of [\:bjr18:]
form a Darboux pair of holomorphic maps by [\:jesz21:, Thm~31],
we also learn that every flat front in hyperbolic space
admits a Weierstrass representation:

\proclaim Corollary.
Any discrete flat front $(X,N):\Sigma_0\to H^3\times S^{2,1}$,
defined as a linear Weingarten net with mixed area Gauss
curvature $K\equiv 1$,
admits a Weierstrass-type representation (\:W:).

%%%%%%%%%%%%%%%%%%%%%%%%%%%%%%%%%%%%%%%%%%%%%%%%%%%%%%%%%%%%%%%
%% references
%%%%%%%%%%%%%%%%%%%%%%%%%%%%%%%%%%%%%%%%%%%%%%%%%%%%%%%%%%%%%%%
\section References
%%%%%%%%%%%%%%%%%%%%%%%%%%%%%%%%%%%%%%%%%%%%%%%%%%%%%%%%%%%%%%%

\bgroup\frenchspacing\parindent=2em
\refitem bl29
 W Blaschke:
 {\it Vorlesungen \"uber Differentialgeometrie III\/};
 Springer Grundlehren XXIX, Berlin (1929)
\refitem bosu08
 A Bobenko, Y Suris:
 {\it Discrete differential geometry. Integrable structure\/};
 Grad Stud Math 98, AMS, Providence (2008)
\refitem br87
 R Bryant:
 {\it Surfaces of mean curvature one in hyperbolic space\/};
 Ast\'erisque 154--155, 321--347 (1987)
\refitem bjr10
 F Burstall, U Hertrich-Jeromin, W Rossman:
 {\it Lie geometry of flat fronts in hyperbolic space\/};
 CR 348, 661--664 (2010)
\refitem bjr12
 F Burstall, U Hertrich-Jeromin, W Rossman:
 {\it Lie geometry of linear Weingarten surfaces\/};
 CR 350, 413--416 (2012)
\refitem bjr18
 F Burstall, U Hertrich-Jeromin, W Rossman:
 {\it Discrete linear Weingarten surfaces\/};
 Nagoya Math J 231, 55--88 (2018)
\refitem bcjpr20
 F Burstall, J Cho, U Hertrich-Jeromin, M Pember, W Rossman:
 {\it Discrete \char 10-nets and Guichard nets\/};
 EPrint arXiv:2008.01447 (2020)
\refitem ce92
 T Cecil:
 {\it Lie sphere geometry\/};
 Springer Universitext, New York (1992)
\refitem cl12
 D Clarke:
 {\it Integrability in submanifold geometry\/};
 PhD thesis, Univ of Bath (2012)
\refitem gmm00
 J Galvez, A Martinez, F Milan:
 {\it Flat surfaces in hyperbolic $3$-space\/};
 Math Ann 316, 419--435 (2000)
\refitem imdg
 U Hertrich-Jeromin:
 {\it Introduction to M\"obius differential geometry\/};
 London Math Soc Lect Note Series 300,
  Cambridge Univ Press, Cambridge (2003)
\refitem jeho17
 U Hertrich-Jeromin, A Honda:
 {\it Minimal Darboux transformations\/};
 Beitr Alg Geom 58, 81--91 (2017)
\refitem jesz21
 U Hertrich-Jeromin, G Szewieczek:
 {\it Discrete cyclic systems and circle congruences\/};
 EPrint arXiv:2104.13441 (2021)
\refitem hrsy12
 T Hoffmann, W Rossman, T Sasaki, M Yoshida:
 {\it Discrete flat surfaces and linear Weingarten surfaces
  in hyperbolic $3$-space\/};
 Trans Amer Math Soc 364, 5605--5644 (2012)
\refitem pe20
 M Pember:
 {\it Weierstrass-type representations\/};
 Geom Dedicata 204, 299--309 (2020)
\refitem ppy21
 M Pember, D Polly, M Yasumoto:
 {\it Discrete Weierstrass-type representations\/};
 EPrint arXiv: 2105.06774 (2021)
\refitem roya16
 W Rossman, M Yasumoto:
 {\it Discrete linear Weingarten surfaces with singularities
  in Riemannian and Lorentzian spaceforms\/};
 Adv Stud Pure Math 78, 383--410 (2018)
\par\egroup
%%%%%%%%%%%%%%%%%%%%%%%%%%%%%%%%%%%%%%%%%%%%%%%%%%%%%%%%%%%%%%%

%%%%%%%%%%%%%%%%%%%%%%%%%%%%%%%%%%%%%%%%%%%%%%%%%%%%%%%%%%%%%%%
%% bottom matter
%%%%%%%%%%%%%%%%%%%%%%%%%%%%%%%%%%%%%%%%%%%%%%%%%%%%%%%%%%%%%%%
\vskip3em\vfill
\bgroup\fn[cmr7]\baselineskip=8pt
\def\addwd{\hsize=.42\hsize}
\def\uhj_gs{\vtop{\addwd
 U Hertrich-Jeromin, G Szewieczek\\
 Institute of Discrete Mathematics and Geometry\\
 Vienna University of Technology\\
 Wiedner Hauptstra\ss{}e 8--10/104\\ A-1040 Vienna (Austria)\\
 Email: udo.hertrich-jeromin@tuwien.ac.at\\
 \phantom{Email: }gudrun@geometrie.tuwien.ac.at
 }}
\def\jdubois{\vtop{\addwd
 J Dubois\\ Equipe ANGE\\ Inria Paris\\
 2, rue Simone Iff\\ CS 42112\\ F-75589 Paris cedex 12\\
 Email: juliette.dubois@inria.fr
 }}
\hbox to \hsize{\hfil \jdubois \hfil \uhj_gs \hfil}\vskip 3ex
\egroup
%%%%%%%%%%%%%%%%%%%%%%%%%%%%%%%%%%%%%%%%%%%%%%%%%%%%%%%%%%%%%%%
\bye